\documentclass[12pt,leqno]{amsart}
\newtheorem{theorem}{Theorem}[section]
\newtheorem{proposition}{Proposition}[section]
\usepackage{amssymb,amsfonts,amsthm,setspace,indentfirst}
\numberwithin{equation}{section}
\usepackage[dvips]{graphics}
\usepackage{epsfig}

\begin{document}
\title[CURVATURE PROPERTIES OF G\"{O}DEL METRIC]{\bf{CURVATURE PROPERTIES OF G\"{O}DEL METRIC}}
\author[R. Deszcz, M. Hotlo\'{s}, J. Je\l owicki, H. Kundu and A. A. Shaikh]{Ryszard Deszcz, Marian Hotlo\'{s}, Jan Je\l owicki, Haradhan Kundu and Absos Ali Shaikh}
\date{}
\address{\noindent\newline Ryszard Deszcz and Jan Je\l owicki, \newline Department of Mathematics,\newline Wroc\l aw University of Environmental and Life Sciences\newline Grunwaldzka 53, 50-357 Wroc\l aw , Poland}
\email{ryszard.deszcz@up.wroc.pl \ jan.jelowicki@up.wroc.pl}
\address{\noindent\newline Marian Hotlo\'{s}\newline Institute of Mathematics and Computer Science\newline Wroc\l aw 
University of Technology\newline  Wybrze\.{z}e Wyspia\'{n}skiego 27, 50-370 Wroc\l aw, Poland}
\email{marian.hotlos@pwr.wroc.pl}
\address{\newline\noindent Haradhan Kundu and Absos Ali Shaikh\newline Department of Mathematics,\newline University of 
Burdwan, Golapbag,\newline Burdwan-713104,\newline West Bengal, India}
\email{kundu.haradhan@gmail.com, \ aask2003@yahoo.co.in, \ aashaikh@mathburuniv.ac.in}
%
%
%
%
\begin{abstract}
The main aim of this article is to investigate the geometric structures admitting by the G\"{o}del spacetime which produces a new class of semi-Riemannian manifolds (see Theorem \ref{thm4.1} and Theorem \ref{thm4.4}). We also consider some extension of G\"{o}del metric (see Example 4.1).
\end{abstract}
\noindent\footnotetext{ 
$\mathbf{2010}$\hspace{5pt}Mathematics\; Subject\; Classification: 53B20, 53B30, 53B50, 53C50, 53C80, 83C57. }
\maketitle
%
\section{Introduction}\label{intro}
\indent In 1949 G\"{o}del \cite{KG} obtained an exact solution of Einstein field equation with a non-zero cosmological constant corresponding 
to a universe in rotation and with an incoherent matter distribution. In that paper he described a metric nowadays called 
G\"{o}del metric as exact and stationary solution of Einstein field equation, which describes a rotating, homogeneous 
but non-isotropic spacetime. Possessing a series of strange properties, it remains still today quite interesting 
mathematically and significant physically. For example, it contains rotating matter but have not singularity, and also 
it is cyclic Ricci parallel \cite{VD}. It is known that the Weyl conformal tensor of the G\"{o}del solution has Petrov type $D$, and G\"{o}del solution is, up to local isometry, the only perfect fluid solution of Einstein field equation admitting five dimensional 
Lie algebra of Killing vectors. G\"{o}del spacetime is geodesically complete, its timelike curves are closed \cite{DM1}. 
Also G\"{o}del spacetime is not globally hyperbolic but diffeomorphic to $\mathbb R^4$ and is simply connected. G\"{o}del metric 
is the Cartesian product of a factor $\mathbb R$ with a three dimensional Lorentzian manifold with signature $(- + + +)$.

\indent G\"{o}del metric and its properties have been studied by various authors to describe 
the G\"{o}del universe. Kundt \cite{WK} studied its geodesics in 1956, and Hawking and Ellis \cite{HE} 
emphasized on coordinates showing its rotational symmetry to draw a nice picture of its dynamics in their book in 1973. 
Malament \cite{DM2} calculated the minimal energy of a closed timelike curve 
of G\"{o}del spacetime. In 2001 Radojevi\'c \cite{RD} 
presented modification of G\"{o}del metric in order to find out some other perfect fluid solutions. Induced matter theory 
and embedding of G\"{o}del universe in five-dimensional Ricci flat space was studied by Fonseca-Neto et. al. \cite{JBF} in 2005. 
Riemann extension of G\"{o}del metric was considered by Dryuma \cite{VD} in 2005, and Dautcourt et. al. \cite{GD} studied light cone 
of G\"{o}del universe. Lanczos spin tensor of G\"{o}del geometry was studied by Garcia-Olivo 
et. al. \cite{GR} in 2006. G\"{o}del metric in various dimensions was studied by G\"{u}rses et. al. \cite{GM}. 
Generalized G\"{o}del metric is given by Plaue et. al. \cite{PS} in 2008.
The G\"{o}del metric \cite{KG} is given by:
\begin{equation}\label{eqn1.1}
ds^2=g_{ij}dx^i dx^j=a^2 \left( -(dx^1)^2 + \frac{1}{2}e^{2x^1}(dx^2)^2 - (dx^3)^2 + (dx^4)^2 + 2e^{x^1}dx^2 dx^4 \right), 
\end{equation}
where $-\infty <x^i< \infty$, $i,j \in \{1,2,3,4 \}$ 
and $a^2 = \frac{1}{2\omega^2}$, $\omega$ is a non-zero real constant, 
which turns out to be the angular velocity, 
as measured by any non-spinning observer located at any one of the dust grains.

\indent The object of the paper is to present the curvature properties of G\"{o}del metric. 
Section 2 deals with semi-Riemannian manifolds with cyclic parallel and Codazzi type Ricci tensor 
and we provide a metric whose Ricci tensor is of Codazzi type but not cyclic parallel (see Example 2.1). However, 
G\"{o}del spacetime is a manifold with cyclic parallel Ricci tensor but the Ricci tensor is not of Codazzi type. 
Section 3 is concerned with rudiments of  pseudosymmetry type manifolds, and 
in the last section we investigate the geometric structures admitting by G\"{o}del metric. Among others, it is shown that G\"{o}del spacetime is neither pseudosymmetric nor Ricci pseudosymmetric but it is quasi-Einstein and a special type of Ricci generalized pseudosymmetric (i.e., $R \cdot R = Q(S,R)$). Although, it is not conformally pseudosymmetric but its Weyl conformal curvature tensor is pseudosymmetric (i.e., $C \cdot C = \frac{\kappa}{6}Q(g,C)$) (see Theorem \ref{thm4.1}). Hence G\"{o}del spacetime induces a new class of semi-Riemannian manifolds which are quasi-Einstein with pseudosymmetric Weyl conformal curvature tensor satisfying $R\cdot R = Q(S,R)$. Finally, we consider some extension of G\"{o}del metric (see Example 4.1).
\section{Manifolds with cyclic parallel and Codazzi type Ricci  tensor}

\indent Let $(M,g)$, $\dim M = n \geq 3$, be a connected paracompact manifold of class $C^{\infty}$ with the metric $g$ 
of signature $(s, n-s)$, $0 \leq s \leq n$. The manifold $(M,g)$ will be called a semi(pseudo)-Riemannian manifold. 
Clearly, if $s = 0$ or $s = n$ then $(M,g)$ is a Riemannian manifold. If $s = 1$ or $s = n-1$ 
then $(M,g)$ is a Lorentzian manifold. Further, let $\nabla$, $R$, $S$ and $\kappa $ be the Levi-Civita connection, 
the curvature tensor, the Ricci tensor and the scalar curvature of the semi-Riemannian manifold $(M,g)$, respectively. 

\indent The semi-Riemannian manifold $(M,g)$ is called locally symmetric if $\nabla R = 0$ (locally $R^{h}_{ijk ,\, l} = 0$), 
which is equivalent to the fact that for each point $x \in M$ the local geodesic symmetry is an isometry. 
For $2$-dimensional manifolds being of locally symmetric and being of constant curvature are equivalent. 
But for $n \geq 3$, the locally symmetric manifolds are a generalization of the manifolds of constant curvature. 
A full classification of locally symmetric manifolds is given by Cartan \cite{CAR} for the Riemannian case, 
and Cahen and Parker (\cite{CAH}, \cite{CAH1}) for the non-Riemannian case. 

\indent The semi-Riemannian manifold $(M,g)$ is said to be Ricci symmetric if $\nabla S = 0$ (locally $S_{ij ,\, k}=0$). 
Every locally symmetric semi-Riemannian manifold is Ricci symmetric but not conversely. However, 
the converse statement is true when $n = 3$. For a compendium of natural symmetries of semi-Riemannian manifolds, 
we refer to \cite{DESZ} and \cite{HaVer02}. 
We mention that Gray in \cite{GRAY}, among other things, investigated various extensions of the class of Ricci symmetric manifolds.
We denote by $\mathcal{A}$  the class of semi-Riemannian manifolds whose Ricci tensor $S$ is cyclic parallel, i.e.,
\begin{equation}\label{eqn-2.1}
(\nabla_X S)(Y,Z) + (\nabla_Y S)(Z,X) + (\nabla_Z S)(X,Y) =  0
\end{equation}
for all vector fields $X$, $Y$, $Z\in \chi(M)$, $\chi(M)$ being the Lie algebra of all smooth vector fields on $M$. 
The local expression of (\ref{eqn-2.1}) is 
$S_{ij, \, k} + S_{jk, \, i} + S_{ki, \, j} = 0$.
A semi-Riemannian manifold satisfying (\ref{eqn-2.1}) is said to be a manifold with cyclic parallel Ricci tensor.
We mention that D'Atri and Nickerson \cite{DAR} proposed to study some class of Riemannian manifolds 
whose curvature tensor satisfies certain conditions of which the first one is equivalent to (\ref{eqn-2.1}). 
 
\indent Evidently, every Ricci symmetric semi-Riemannian manifold is a manifold with cyclic parallel Ricci tensor but not conversely. 
However, the converse statement is true if the Ricci tensor is a Codazzi tensor. We recall that an $(0,2)$-symmetric tensor 
$B$ 
is said to be a Codazzi tensor if it satisfies the Codazzi equation, i.e. 
$(\nabla_X B)(Y,Z) = (\nabla_Y B)(X,Z)$. 
The local expression of the last equation is 
$B_{ij, \, k} = B_{kj, \, i}$, 
where $B_{ij}$ are the local components of the tensor $B$.
A Codazzi tensor is trivial if it is a constant multiple of the metric tensor \cite{DER2}.
We denote by $\mathcal{B}$
the class of semi-Riemannian manifolds with Ricci tensor $S$ as Codazzi tensor, i.e. $S$ 
satisfies
$(\nabla_XS)(Y,Z) = (\nabla_YS)(X,Z)$,
for all vector fields $X$, $Y$, $Z\in \chi(M)$. Every Ricci symmetric semi-Riemannian manifold 
is of class $\mathcal{B}$ but not conversely. Codazzi tensors are of great interest in the geometric literature 
and have been studied by several authors, as Berger and Ebin \cite{BER}, Bourguignon \cite{BOUR}, Derdzinski (\cite{DER}, \cite{DER1}),
Derdzinski and Shen \cite{DER2}, Ferus \cite{FER}, Simon \cite{SIM}; a compendium of results is reported in Besse's book \cite{BEESE}.

\indent We note that every semi-Riemannian manifold of constant curvature and hence Einstein semi-Riemannian manifold is of class $\mathcal{A}$ 
as well as of $\mathcal{B}$. 
We note that the scalar curvature $\kappa $ of every semi-Riemannian manifold of the class $\mathcal{A}$ or $\mathcal{B}$
is constant.
We note that G\"{o}del spacetime is of class $\mathcal A$ but not of class $\mathcal B$.

\indent It is known that Cartan hypersurfaces are Riemannian manifolds, with non-parallel Ricci tensor, 
satisfying the generalized Einstein metric condition 
(\ref{eqn-2.1}) (\cite{Naka}, Theorem 4.1). As it was noted in \cite{D5}(p. 109), the Cartan hypersurfaces do not satisfy 
\begin{equation}\label{1}
\nabla_Z\left(S(X, Y)-\frac{\kappa }{2(n-1)}g(X, Y)\right)=\nabla_Y\left(S(X, Z)-\frac{\kappa }{2(n-1)}g(X, Z)\right).
\end{equation}
We mention that (\ref{1}) is presented in the Table 1, pp. 432-433 of \cite{BEESE}. We also refer to  \cite{BEESE} for results
on Riemannian manifolds satisfying (\ref{1}).
This was  noted that 
Codazzi tensors occur naturally in the study of harmonic Riemannian manifolds. The Ricci tensor is a Codazzi tensor 
if and only if $div\, R=0$ i.e., if and only if the manifold has harmonic curvature tensor \cite{BEESE}.

\indent We note that in a $3$-dimensional Riemannian manifold $(M,g)$, the following conditions:\\
(a) $(M,g)$ is locally symmetric, 
(b) $(M,g)$ is Ricci symmetric and 
(c) $(M,g)$ is a conformally flat manifold with cyclic parallel Ricci tensor,
are equivalent \cite{SB}.
Also for a Riemannian manifold $(M,g)$, $n \geq 4$, the following conditions:
(a) $(M,g)$ is Ricci symmetric, 
(b) $(M,g)$ is a manifold with cyclic parallel Ricci tensor and harmonic conformal curvature tensor, and
(c) $(M,g)$ is a manifold with cyclic parallel and Codazzi type Ricci tensor, are equivalent \cite{SB}.
\newline

\noindent
{\bf{Example 2.1}}
(i) Let $M = \{(x^1,x^2,x^3,x^4) : x^i>0,i = 1, 2, 3, 4\}$ be the subset of 
$\mathbb R^4$ endowed with the metric $g$ defined by
$ds^2= \varepsilon (dx^1)^2+x^1\left(  (dx^2)^2 +  (dx^3)^2 + (dx^4)^2 \right)$, $\varepsilon = \pm 1$.
It is easy to check that $(M,g)$ is a conformally flat 
quasi-Einstein Riemannian manifold,
$\mbox{rank}\, (S - \frac{1}{4 (x^{1})^{2}}\, g) = 1$,
whose scalar curvature $\kappa$ is equal to zero, 
and the 
Ricci tensor $S$ is of Codazzi type but not cyclic parallel.
Moreover, we have 
$R \cdot R =  Q(S,R)$ 
and
$R \cdot R = L\, Q(g,R)$, $L = - \frac{1}{4 \varepsilon (x^{1})^{2}}$.
\newline
(ii) Let $M = \{(x^1,x^2, \ldots ,x^5) : x^i>0,i = 1, 2, \ldots , 5 \}$ be the subset of 
$\mathbb R^5$ endowed with the metric $g$ defined by
$ds^2= \varepsilon (dx^1)^2+x^1\left( (dx^2)^2 + \ldots  + (dx^5)^2 \right)$, $\varepsilon = \pm 1$.
It is easy to check that $(M,g)$ is a conformally flat 
quasi-Einstein Riemannian manifold,
$\mbox{rank}\, (S - \frac{1}{2 (x^{1})^{2}}\, g) = 1$, 
whose scalar curvature $\kappa$ is non-zero, $\kappa = \frac{1}{ \varepsilon (x^{1})^{2}}$,	
and the Ricci tensor $S$ is not of Codazzi type and 
not cyclic parallel.
Moreover, we have
$R \cdot R =  Q(S,R)$ 
and 
$R \cdot R = L\, Q(g,R)$, $L = - \frac{1}{4 \varepsilon (x^{1})^{2}}$.
\section{Pseudosymmetry type curvature conditions}

\indent We define on a semi-Riemannian manifold $(M,g)$, $n \geq 3$, the endomorphisms $X\wedge_A Y$,
$\mathcal{R}(X,Y)$, $\mathcal{C}(X,Y)$, $\mathcal{K}(X,Y)$ and $conh(\mathcal{R})$
by (\cite{D0}, \cite{DGHSaw}, 
\cite{GLOG2},
\cite{Ishii},
\cite{Kow2})
\begin{eqnarray*}
(X\wedge_A Y)Z &=& A(Y,Z)X-A(X,Z)Y,\\
\mathcal{R}(X,Y)Z &=& [\nabla_X,\nabla_Y]Z-\nabla_{[X,Y]}Z,\\
\mathcal{C}(X,Y) &=& 
\mathcal{R}(X,Y)-\frac{1}{n-2}(X\wedge_g \mathcal{L} Y + \mathcal{L} X\wedge_g Y - \frac{\kappa }{n-1}X\wedge_g Y),\\
\mathcal{K}(X,Y) &=& \mathcal{R}(X,Y) - \frac{\kappa }{n (n-1)} X\wedge_g Y,\\
{conh(\mathcal{R})}(X,Y) &=& {\mathcal{R}}(X,Y) 
- \frac{1}{n-2}(X \wedge _{g} {\mathcal{L}}Y + {\mathcal{L}}X \wedge _{g} Y ),
\end{eqnarray*}
respectively, where $A$ is an $(0,2)$-tensor on $M$, $X,Y,Z\in\chi(M)$. The Ricci operator $\mathcal{L}$ is defined
by $g(X,\mathcal{L} Y) = S(X,Y)$, where $S$ is the Ricci tensor and $\kappa $
the scalar curvature of $(M,g)$, respectively. 
The tensor $S^{2}$ is defined by $S^{2}(X,Y) = S(X,\mathcal{L} Y)$.
Further, we define the Gaussian curvature tensor $G$,
the Riemann-Christoffel curvature tensor $R$, the Weyl conformal
curvature tensor $C$, concircular curvature tensor $K$ and conharmonic curvature tensor $conh(R)$ of $(M,g)$, 
by (\cite{D0}, 
\cite{DGHSaw}, 
\cite{GLOG2}, \cite{Ishii})
\begin{eqnarray*}
G(X_1,X_2,X_3,X_4)     &=& g((X_1\wedge_g X_2)X_3,X_4),\\
R(X_1,X_2,X_3,X_4)     &=& g(\mathcal{R}(X_1,X_2)X_3,X_4),\\
C(X_1,X_2,X_3,X_4)     &=& g(\mathcal{C}(X_1,X_2)X_3,X_4),\\
K(X_1,X_2,X_3,X_4)     &=& g(\mathcal{K}(X_1,X_2)X_3,X_4),\\
conh(R)(X_1,X_2,X_3,X_4) &=& g(conh(\mathcal{R})(X_{1},X_{2})X_{3},X_{4}) ,
\end{eqnarray*}
respectively. 
For $(0,2)$-tensors $A$ and $B$ we define their Kulkarni-Nomizu product $A\wedge B$
by (see e.g. \cite{DGHSaw}, 
\cite{GLOG2}, 
\cite{Kow2})
\begin{eqnarray*}
(A\wedge B)(X_1,X_2;X,Y)&=&A(X_1,Y)B(X_2,X)+A(X_2,X)B(X_1,Y)\\
&-&A(X_1,X)B(X_2,Y)-A(X_2,Y)B(X_1,X).
\end{eqnarray*}
We note that the Weyl conformal curvature tensor $C$ can be presented in the following form
\begin{equation}
\label{eqn2.1}
C = R-\frac{1}{n-2}g\wedge S + \frac{\kappa }{(n-2)(n-1)}G.
\end{equation}
For an $(0,k)$-tensor $T$, $k\geq 1$ and a symmetric $(0,2)$-tensor $A$ we define the
$(0,k)$-tensor $A \cdot T$ and the $(0,k+2)$-tensors $R \cdot T$, $C \cdot T$ and $Q(A, T)$ by
\begin{eqnarray*}
&&(A\cdot T)(X_1,\cdots, X_k)=-T(\mathcal{A}X_1,X_2,\cdots,X_k)-\cdots-T(X_1,X_2,\cdots,\mathcal{A}X_k),\\
&&(R\cdot T)(X_1,\cdots, X_k;X,Y)=(\mathcal{R}(X,Y)\cdot T)(X_1,\cdots, X_k)\\
&&=-T(\mathcal{R}(X,Y)X_1,X_2,\cdots,X_k)-\cdots-T(X_1,\cdots,X_{k-1},\mathcal{R}(X,Y)X_k),\\
&&(C\cdot T)(X_1,\cdots, X_k;X,Y)=(\mathcal{C}(X,Y)\cdot T)(X_1,\cdots, X_k)\\
&&=-T(\mathcal{C}(X,Y)X_1,X_2,\cdots,X_k)-\cdots-T(X_1,\cdots,X_{k-1},\mathcal{C}(X,Y)X_k),\\
&&Q(A, T)(X_1,\cdots, X_k;X,Y)=((X\wedge_A Y)\cdot T)(X_1,\cdots, X_k)\\
&&=-T((X\wedge_A Y)X_1,X_2,\cdots,X_k)-\cdots-T(X_1,\cdots,X_{k-1},(X\wedge_A Y)X_k),
\end{eqnarray*}
where $\mathcal{A}$ is the endomorphism of $\chi(M)$ defined by
$g(\mathcal{A}X,Y) = A(X,Y)$. Putting in the above formulas $T = R, T = S$, $T = C$ or $T = K$, $A = g$ or $A = S$, 
we obtain the tensors: $R\cdot R$, $R\cdot S$, $R\cdot C$, $R\cdot K$, $C\cdot R$, $C\cdot S$, $C\cdot C$, 
$C\cdot K$, $Q(g, R)$, $Q(g, S)$, $Q(g, C)$, $Q(g, K)$, $Q(S, R)$, $Q(S, C)$, $Q(g, K)$, 
$S\cdot R$, $S\cdot C$, $S\cdot K$, 
$conh(R) \cdot conh(R)$, $conh(R) \cdot R$, $R \cdot conh(R)$, $conh(R) \cdot S$, 
etc.
The tensor $Q(A,T)$ is called the Tachibana tensor of the tensors $A$ and $T$, or the Tachibana
tensor for short (\cite{DGPSS}). We like to point out that in some papers, $Q(g,R)$ is called the
Tachibana tensor 
(see e.g. \cite{DGPSS}, \cite{HaVer01}, \cite{JaHaSenVer},  \cite{JaHaPetroVer}).
We also have

\begin{proposition} (cf. \cite{DGH01}) For any semi-Riemannian manifold $(M,g)$, $n \geq 4$, we have 
\begin{eqnarray}
conh(R) \cdot S &=& C \cdot S - \frac{ \kappa }{(n-2)(n-1)}\, Q( g, S),\nonumber\\ 
R \cdot conh(R) &=&  R \cdot C , \nonumber\\
conh(R) \cdot R 
&=&
C \cdot R - \frac{ \kappa }{(n-2)(n-1)}\, Q(g,R) ,\nonumber\\ 
conh(R) \cdot conh(R) 
&=&  
C \cdot C - \frac{ \kappa }{(n-2)(n-1)}\, Q(g, C) .
\label{conh02}
\end{eqnarray}
\end{proposition}

\indent A semi-Riemannian manifold $(M, g)$, $n \geq 3$, satisfying the condition
\begin{equation}
\label{4.3.001}
R \cdot R = 0 
\end{equation}
is called semisymmetric (\cite{sz}). 
We mention that 
non-conformally flat and non-locally symmetric semi-Riemannian manifolds
having parallel Weyl conformal curvature tensor are semisymmetric (\cite{DerRo01}, Theorem 9), 
their scalar curvature is equal to zero (\cite{DerRo01}, Theorem 7) 
and the Ricci tensor is a Codazzi tensor (\cite{DerRo01}, eq. (6)). 
We refer to \cite{DerRo03}-\cite{DerRo06} for the recent results on semi-Riemannian manifolds with parallel Weyl conformal curvature tensor
and, in particular, for classification results. 
Semi-Riemannian warped products having parallel Weyl conformal curvature tensor were investigated in \cite{Hotlos}. 
We also mention that, recently, conformally semisymmetric manifolds and special semisymmetric Weyl conformal tensors are studied in \cite{es}.
Another important subclass of semisymmetric semi-Riemannian manifolds form manifolds satisfying
\begin{equation}
\label{4.3.00177}
\nabla \nabla R = 0 .
\end{equation}
We refer to \cite{BSS} and \cite{Senovilla} and references therein for results on manifolds satisfying (\ref{4.3.00177}).

\indent A semi-Riemannian manifold $(M,g)$, $n \geq 3$, is said to be pseudosymmetric \cite{DEG} 
if the tensor $R \cdot R$ and the Tachibana tensor $Q(g,R)$ are linearly dependent at every point of $M$. 
This is equivalent to 
\begin{equation}
\label{4.3.002}
R \cdot R = L_R\, Q(g,R) 
\end{equation}
on $U_R = \{ x \in M : R - \frac{\kappa }{n(n-1)}G \neq 0\ \mbox{at}\ x \}$, 
where $L_R$ is some function on this set.  
We refer to
\cite{P92},
\cite{DES2}, \cite{DESZ}, \cite{HaVer01} and \cite{HaVer02} for  surveys 
on such manifolds. 
In particular, a geometrical interpretation of  pseudosymmetric manifolds, 
in the Riemannian case, 
is given in \cite{HaVer01}.

\indent We note that \cite{DEG} is the first publication, in which  
a semi-Riemannian manifold satisfying (\ref{4.3.002}) 
was named the pseudosymmetric manifold. 
In \cite{DEG} pseudosymmetric warped products with $1$-dimensional base manifold and $(n-1)$-dimensional fibre, $n \geq 4$, 
which is not a semi-Riemannian space of constant curvature, were investigated. 
In \cite{D7} it was shown that
hypersurfaces in spaces of constant curvature, with exactly two distinct principal curvatures at every point, 
are pseudosymmetric.
Thus in particular, Cartan's and Schouten's investigations of quasi-umbilical
hypersurfaces in spaces of constant curvature are closely related to pseudosymmetric
manifolds (cf. \cite{DESZ}). 
It is clear that every semisymmetric manifold is pseudosymmetric.
However, the converse statement is not true.
For instance, the Schwarzschild spacetime, the Kottler spacetime
and the Reissner-Nordstr\"{o}m spacetime satisfy (\ref{4.3.002}) with non-zero function $L_{R}$ \cite{DeVerVra} 
(see also \cite{P109}, \cite{HaVer03}). 
We also mention that 
Friedmann-Lema{\^{\i}}tre-Robertson-Walker spacetimes are 
pseudosymmetric (cf. \cite{DESZ}).
It is well-known that the
Schwarzschild spacetime was discovered in 1916 by Schwarzschild, during
his study on solutions of Einstein's equations. It seems that the Schwarzschild spacetime 
is the first example of a non-semisymmetric, pseudosymmetric warped product. 
Finally, we note that
(\ref{4.3.002}) is equivalent to (e.g. see \cite{DES2})
\begin{equation}
\label{4.3.004}(R - L_{R}\, G) \cdot (R  - L_{R}\, G) = 0 .
\end{equation}

\indent We also note that in \cite{CHA} Chaki introduced another kind of pseudosymmetry. However,
both notions of pseudosymmetry are not equivalent. 
Throughout the paper we will confine the pseudosymmetry related to (\ref{4.3.002}).

\indent A semi-Riemannian manifold $(M,g)$, $n \geq 3$, is said to be Ricci-pseudosymmetric (\cite{DES1}, \cite{DH}) 
if the tensor $R \cdot S$ and the Tachibana tensor $Q(g,S)$ are linearly dependent at every point of $M$. 
Thus the manifold $(M,g)$ is Ricci-pseudosymmetric if and only if 
\begin{equation}
\label{4.3.1}
R \cdot S = L_{S}\, Q(g,S)
\end{equation}
holds on $U_{S} = \{ x \in M: S - \frac{r}{n}\, g \neq 0\ \mbox{at}\ x \}$, 
where $L_{S}$ is  some function on this set. We note that $U_{S} \subset U_{R}$. 
It is easy to check that (\ref{4.3.1}) is equivalent to 
\begin{equation}
\label{4.3.006}
(R - L_{S}\, G) \cdot (S  - L_{S}\, g) = 0 .
\end{equation}
We refer to \cite{P92}, \cite{DGHSaw}, \cite{DESZ}, \cite{D5} and \cite{HaVer01}  
for  surveys and comments on such manifolds. A geometrical interpretation of Ricci-pseudosymmetric manifolds, 
in the Riemannian case, is given in \cite{JaHaSenVer}. 
It is clear that every pseudosymmetric semi-Riemannian manifold is Ricci-pseudosymmetric.
However, the converse statement is not true. For instance,
the Cartan hypersurfaces of dimension $6$, $12$ or $24$
are non-quasi-Einstein and non-pseudosymmetric Ricci-pseudosymmetric manifolds 
(\cite{DY01}, see also \cite{DGHSaw}, \cite{D5}). 
The $3$-dimensional Cartan hypersurface 
is a quasi-Einstein pseudosymmetric manifold \cite{D7}. 
We mention that recently quasi-Einstein 
Ricci-pseudosymmetric hypersurfaces in semi-Riemannian spaces of constant curvature 
were investigated in \cite{DHS01}.

\indent A semi-Riemannian manifold $(M,g)$, $n \geq 4$, is said to be conformally pseudosymmetric (\cite{DES2}, \cite{D30})
if the tensor $R \cdot C$ and the Tachibana tensor $Q(g,C)$ are linearly dependent at every point of $M$. 
Again a semi-Riemannian manifold $(M,g)$, $n \geq 4$, is said to be a manifold with pseudosymmetric Weyl conformal curvature tensor
(\cite{D30}, \cite{DES2})
if the tensor $C \cdot C$ and the Tachibana tensor $Q(g,C)$ are linearly dependent at every point of $M$. 
This is equivalent to 
\begin{equation}
\label{4.3.012}
C \cdot C = L_{C}\, Q(g,C)
\end{equation}
on $U_C = \{ x \in M : C \neq 0\ \mbox{at}\ x \}$, 
where $L_{C}$ is some function on this set. We note that $U_{C} \subset U_{R}$. 
It is easy to check that (\ref{4.3.012}) is equivalent to 
\begin{equation}
\label{4.3.0010}
(C - L_{C}\, G) \cdot (C  - L_{C}\, G) = 0 .
\end{equation}
Using (\ref{eqn2.1}), we also can check 
that 
(\ref{4.3.004}) and (\ref{4.3.0010}) are equivalent on every Einstein manifold. 
As it was stated in \cite{D30}, any warped product $M_{1} \times_{F} M_{2}$, 
with $\dim\, M_{1} = \dim M_{2} = 2$, satisfies (\ref{4.3.012}). Thus in particular,
the Schwarzschild spacetime, the Kottler spacetime
and the Reissner-Nordstr\"{o}m spacetime satisfy (\ref{4.3.012}).

\indent A semi-Riemannian manifold $(M, g)$, $n \geq 3$, is said to be Ricci-generalized pseudosymmetric 
(\cite{DEF}, \cite{DEF1}) 
if at every point of $M$ the tensor $R\cdot R$ and the Tachibana tensor $Q(S,R)$ are linearly dependent. 
Hence $(M, g)$ is Ricci-generalized pseudosymmetric if and only if 
\begin{equation}
\label{4.3.00777}
R\cdot R=L\, Q(S,R)
\end{equation}
holds on  $U=\{x\in M : Q(S,R) \neq 0\ \mbox{at}\ x \}$, where $L$ is some function on this set. 
An important subclass of Ricci-generalized pseudosymmetric manifolds is formed by the manifolds realizing 
the condition (\cite{DEF}, \cite{DG})
\begin{equation}
\label{4.3.00444}
R \cdot R=Q(S,R) . 
\end{equation}

\indent At the end of this section we also present some other curvature conditions closely related to the above presented conditions.
Namely, it was stated in \cite{32}, on every hypersurface $M$ immersed isometrically 
in a semi-Riemannian space of constant curvature $N$, $\dim\, N = n+1$, $n \geq 4$, we have
\begin{equation}
\label{4.3.00555}
R \cdot R - Q(S,R) = - \frac{(n-2) \widetilde{\kappa }}{n(n+1)}\, Q(g,C) , 
\end{equation}
where $\widetilde{\kappa }$ is the scalar curvature of the ambient space. 
It is clear, that if the ambient space is a semi-Euclidean space then (\ref{4.3.00555}) 
reduces to (\ref{4.3.00444}). We also note that any warped product $M_{1} \times_{F} M_{2}$, 
with $\dim\, M_{1} = 1$, $\dim M_{2} = 3$, satisfies \cite{49}
\begin{equation}
\label{4.3.00666}
R \cdot R - Q(S,R) = L\, Q(g,C) , 
\end{equation}
for some function $L$ on $U_{C}$. Thus generalized Robertson-Walker spacetimes fulfills  
(\ref{4.3.00666}). In particular, Friedmann-Lema{\^{\i}}tre-Robertson-Walker spacetimes satisfy (\ref{4.3.00444}).
We mention that the Vaidya spacetime also satisfies (\ref{4.3.00666}) (\cite{Kow2}, Example 5.2).

\indent The conditions: (\ref{4.3.002}), (\ref{4.3.1}), (\ref{4.3.012}), (\ref{4.3.00777}) and (\ref{4.3.00666}) 
or other conditions of this kind are called conditions of pseudosymmetry type. 
We refer to \cite{P92},  \cite{D0}, \cite{DES2}, \cite{DGHSaw} and \cite{DESZ} 
for surveys on semi-Riemannian manifolds satisfying such conditions. In particular, 
we refer to \cite{CH-DDGP} for recent results on quasi-Einstein manifolds 
satisfying curvature conditions of this kind.

\indent It is easy to check that on 
every Ricci-pseudosymmetric manifold 
the following condition is satisfied 
(\cite{ACDE1998}, Lemma 3.3; \cite{DGHS}, Proposition 3.1(iv))
\begin{equation}
\label{A09}
R ( {\mathcal{L}}X, Y,Z,W) 
+
R ( {\mathcal{L}}Z, Y,W,X)
+
R ( {\mathcal{L}}W, Y,X,Z) = 0 . 
\end{equation}
Semi-Riemannian manifolds $(M,g)$, $n \geq 3$, satisfying 
(\ref{A09}) are called Riemann compatible (\cite{23}, \cite{22}).
We also note that (\ref{A09}) remains invariant under geodesic mappings.  
In \cite{DGPSS} (Proposition 2.1) it was proved that 
(\ref{A09}) holds on every manifold satisfying 
(\ref{4.3.00666}). Thus in particular, manifolds satisfying 
(\ref{4.3.00444}) are also Riemann compatible 
(\cite{DEF}, Lemma 2.2 (i)).
We refer to \cite{DGJP-TZ}, \cite{24} and \cite{21} for further results on Riemann compatible manifolds.

\section{G\"{o}del metric admitting geometric structures}\label{Godel}

Let on $\mathbb R^4$ be given the G\"{o}del metric $g$ defined by (\ref{eqn1.1}) and let 
$h,i,j,k,l.m \in \{ 1,2,3,4 \}$. 
From (\ref{eqn1.1}) the non-zero components $\Gamma ^{h}_{ij}$ of the Christoffel symbols of second kind 
of $g$ are given by \cite{KG}:
$$\Gamma^{4}_{12}=\frac{e^{x^1}}{2},\,\,\, \Gamma^{2}_{14}=-e^{-x^1}, \,\,\, 
\Gamma^{4}_{14}=1, \,\,\, \Gamma^{1}_{22}=\frac{e^{2x^1}}{2}, \,\,\, \Gamma^{1}_{24}=\frac{e^{x^1}}{2}.$$
Further, the non-zero components 
$R_{hijk}$ and $S_{ij}$ 
of the Riemann-Christoffel curvature tensor $R$ and the Ricci tensor $S$, respectively, and the scalar curvature $\kappa$ 
are given by \cite{KG}:
$$R_{1212}=\frac{3}{4}a^2e^{2x^1}, \,\,\, R_{1214}=\frac{1}{2}a^2e^{x^1}, \,\,\, 
R_{1414}=\frac{a^2}{2}, \,\,\, R_{2424}=\frac{1}{4}a^2e^{2x^1},$$ 
$$S_{22}=e^{2x^1}, \,\,\, S_{24}=e^{x^1}, \,\,\, S_{44}=1 \mbox{ and } \kappa =\frac{1}{a^2}.$$
Again the non-zero components $R_{hijk, l}$ and $S_{ij, l}$ of the covariant derivatives of the Riemann-Christoffel curvature tensor
$\nabla R$ and the Ricci tensor $\nabla S$, respectively, are given by:
$$R_{1212, 1}=a^2e^{2x^1}, \,\,\, R_{1214, 1}=\frac{1}{2}a^2e^{x^1}, \,\,\, 
R_{1224, 2}=\frac{1}{4}a^2e^{3x^1},$$ 
$$S_{12, 2}=-\frac{e^{2x^1}}{2}, \,\,\, S_{14, 2}=-\frac{e^{x^1}}{2}, \,\,\, 
S_{22, 1}=e^{2x^1}, \,\,\, S_{24, 1}=\frac{e^{x^1}}{2}.$$
The non-zero components $C_{hijk}$ of the Weyl conformal curvature tensor $C$ are given below:
$$C_{1212}=\frac{1}{3}a^2e^{2x^1}, \,\,\, C_{1214}=\frac{1}{6}a^2e^{x^1}, \,\,\, 
C_{1414}=\frac{a^2}{6}, \,\,\, C_{2424}=\frac{1}{12}a^2e^{2x^1},$$
$$C_{1313}=-\frac{1}{6}a^2, \,\,\, C_{2323}=-\frac{5}{12}a^2e^{2x^1}, \,\,\, 
C_{2334}=\frac{1}{3}a^2e^{x^1}, \,\,\, C_{3434}=-\frac{1}{3}a^2.$$
The non-zero components $(R \cdot R)_{hijklm}$ of the tensor $R \cdot R$ are given below:
$$(R\cdot R)_{122412}=\frac{1}{4}a^2e^{3x^1}, \,\,\,   (R \cdot R)_{122414}=\frac{1}{4}a^2e^{2x^1},$$
$$(R \cdot R)_{121224}=-\frac{1}{2}a^2e^{3x^1}, \,\,\, (R \cdot R)_{121424}=-\frac{1}{4}a^2e^{2x^1}.$$
The non-zero components $Q(S,R)_{hijklm}$ of the Tachibana tensor $Q(S, R)$ are given below:
$$Q(S, R)_{122412}=\frac{1}{4}a^2e^{3x^1}, \,\,\,  Q(S, R)_{122414}=\frac{1}{4}a^2e^{2x^1},$$
$$Q(S, R)_{121224}=-\frac{1}{2}a^2e^{3x^1}, \,\,\, Q(S, R)_{121424}=-\frac{1}{4}a^2e^{2x^1}.$$
The non-zero components $(C \cdot C)_{hijklm}$ of the tensor $C \cdot C$ are given below:
\begin{eqnarray*}
& & 
(C \cdot C)_{122412}=\frac{1}{24}a^2e^{3x^1}, \,\,\,
(C \cdot C)_{132312}=\frac{1}{12}a^2e^{2x^1}, \,\,\, 
(C \cdot C)_{133412}=-\frac{1}{12}a^2e^{x^1},\\ 
& & 
(C \cdot C)_{122313}=-\frac{1}{8}a^2e^{2x^1}, \,\,\, 
(C \cdot C)_{123413}=\frac{1}{12}a^2e^{x^1}, \,\,\, 
(C \cdot C)_{142313}=-\frac{1}{12}a^2e^{x^1},\\
& & 
(C \cdot C)_{143413}=\frac{1}{12}a^2, \,\,\,
(C \cdot C)_{122414}=\frac{1}{24}a^2e^{2x^1}, \,\,\, 
(C \cdot C)_{132314}=\frac{1}{12}a^2e^{x^1},\\
& & 
(C \cdot C)_{133414}=-\frac{1}{12}a^2, \,\,\,
(C \cdot C)_{121323}=\frac{1}{24}a^2e^{2x^1}, \,\,\, 
(C \cdot C)_{232423}=-\frac{1}{24}a^2e^{3x^1},\\
& & 
(C \cdot C)_{243423}=\frac{1}{24}a^2e^{2x^1}, \,\,\,
(C \cdot C)_{121224}=-\frac{1}{12}a^2e^{3x^1}, \,\,\, 
(C \cdot C)_{121424}=-\frac{1}{24}a^2e^{2x^1},\\
& & 
(C \cdot C)_{232324}=\frac{1}{12}a^2e^{3x^1}, \,\,\,
(C \cdot C)_{233424}=-\frac{1}{24}a^2e^{2x^1}.
\end{eqnarray*}
The non-zero components $Q(g,C)_{hijklm}$ of the Tachibana tensor $Q(g, C)$ are given below:
\begin{eqnarray*}
& & 
Q(g, C)_{122412}=\frac{1}{4}a^4e^{3x^1}, \,\,\, Q(g, C)_{132312}=\frac{1}{2}a^4e^{2x^1}, \,\,\, Q(g, C)_{133412}=-\frac{1}{2}a^4e^{x^1},\\
& &  
Q(g, C)_{122313}=-\frac{3}{4}a^4e^{2x^1}, \,\,\, Q(g, C)_{123413}=\frac{1}{2}a^4e^{x^1}, \,\,\, Q(g, C)_{142313}=-\frac{1}{2}a^4e^{x^1},\\
& &  
Q(g, C)_{143413}=\frac{1}{2}a^4, \,\,\, Q(g, C)_{122414}=\frac{1}{4}a^4e^{2x^1}, \,\,\, Q(g, C)_{132314}=\frac{1}{2}a^4e^{x^1},\\
& & 
Q(g, C)_{133414}=-\frac{1}{2}a^4, \,\,\, Q(g, C)_{121323}=\frac{1}{4}a^4e^{2x^1}, \,\,\, Q(g, C)_{232423}=-\frac{1}{4}a^4e^{3x^1},\\
& & 
Q(g, C)_{243423}=\frac{1}{4}a^4e^{2x^1}, \,\,\, Q(g, C)_{121224}=-\frac{1}{2}a^4e^{3x^1}, \,\,\, Q(g, C)_{121424}=-\frac{1}{4}a^4e^{2x^1},\\
& & 
Q(g, C)_{232324}=\frac{1}{2}a^4e^{3x^1}, \,\,\, Q(g, C)_{233424}=-\frac{1}{4}a^4e^{2x^1}.
\end{eqnarray*}
Thus we see that $\mathbb{R}^4$ equipped with the G\"{o}del metric $g$ has the following curvature properties:\\
(i) The Ricci tensor is cyclic parallel \cite{VD}, the rank of the Ricci tensor $S$ is $1$ \cite{KG}, precisely,
\begin{equation}
\label{aaa01}
S = \kappa\, \omega \otimes \omega, \,\,\,  \kappa = \frac{1}{a^2}, \,\,\, \omega 
= ( \omega _{1}, \omega _{2}, \omega _{3}, \omega _{4})
= (0, ae^{x^1}, 0, a), 
\end{equation}
and
the vector field $X$ corresponding to $1$-form $\omega $ is given by $X=(0, 0, 0, \frac{1}{a})$,\\
(ii) $R \cdot R = Q(S,R)$,\\
(iii) $C \cdot C = \frac{\kappa }{6} \, Q(g,C)$,\\
(iv) $3 R \cdot K - 2\, Q(S, K) = Q(S,C)$.\\
G\"{o}del metric also realizes the following pseudosymmetric type conditions:\\
(v) \begin{eqnarray*}
(2 a^2 L_1 + \frac{2}{3} L_2)(R\cdot C + C\cdot R) 
&=& (-\frac{2}{3} L_1 + \frac{1}{9 a^2} L_2)\Big(Q(g,R) - 3 Q(S,R)\Big)\\
&+& L_1 Q(g,C) + L_2 Q(S,C),
\end{eqnarray*}
where $L_1$ and $L_2$ are some functions. This condition implies that $R\cdot C$, $C\cdot R$, $Q(g,R)$, $Q(S,R)$, $Q(g,C)$ and $Q(S,C)$ are linearly dependent.\\
(vi) \begin{eqnarray*}
(L_1 + L_2)(C\cdot K + K\cdot C) &=& (\frac{1}{12 a^2} L_1 + \frac{7}{12 a^2} L_2)Q(g,C) + L_1 Q(S,C)\\
&-& \frac{1}{2 a^2} L_2 Q(g,K) +  L_2 Q(S,K),
\end{eqnarray*}
where $L_1$ and $L_2$ are some functions. This condition implies that $C\cdot K$, $K\cdot C$, $Q(g,C)$, $Q(S,C)$, $Q(g,K)$ and $Q(S,K)$ are linearly dependent.\\
(vii) \begin{eqnarray*}
&&(-\frac{2}{5} L_1 + \frac{12 a^2}{5} L_2)(K\cdot conh(R) + conh(R)\cdot K)\\
&=& (\frac{1}{30 a^2} L_1 - \frac{6}{5} L_2)Q(g,K) + L_1 Q(S,K) + L_2 Q(g,conh(R))\\
&+& (-\frac{7}{5} L_1 + \frac{12 a^2}{5} L_2) Q(S,conh(R)),
\end{eqnarray*}
where $L_1$ and $L_2$ are some functions. This condition implies that $K\cdot conh(R)$, $conh(R)\cdot K$, $Q(g,K)$, $Q(S,K)$, $Q(g,conh(R))$ and $Q(S,conh(R))$ are linearly dependent.\\
(viii) \begin{eqnarray*}
2 a^2 L_1 (R\cdot conh(R) + conh(R)\cdot R) &=& -L_1 \Big(Q(g,R)-Q(g,conh(R))\Big)\\
&+& L_2 Q(S,R) + (2 a^2 L_1 - L_2) Q(S,conh(R)),
\end{eqnarray*}
where $L_1$ and $L_2$ are some functions. This condition implies that $R\cdot conh(R)$, $conh(R)\cdot R$, $Q(g,R)$, $Q(S,R)$, $Q(g,conh(R))$ and $Q(S,conh(R))$ are linearly dependent.\\
\indent We note that the condition $\mbox{rank}\, S = 1$ holds at a point of a semi-Riemannian manifold $(M,g)$, $n \geq 3$, if and only if 
$S\wedge S = 0$ at this point. Further, it is easy to check that (iv) is an immediate consequence of (ii) and $S\wedge S = 0$ 
and the definitions of the tensors $C$ and $K$.  

\indent We also note that the condition 
(\ref{4.3.00444})
holds at every point $x$ of a semi-Riemannian manifold $(M,g)$, $n \geq 3$, at which the condition
$$\omega(X_1)\mathcal{R}(X_2,X_3)+\omega(X_2)\mathcal{R}(X_3,X_1)+\omega(X_3)\mathcal{R}(X_1,X_2)=0$$
is satisfied, where $\omega$ is a non-zero covector at $x$ (\cite{DG}, Theorem 3.1, \cite{DEF}, p. 110). 
However, in case of G\"{o}del metric
(\ref{4.3.00444})
holds, but does not satisfy the above condition. 
Since the G\"{o}del metric is a product metric of a $3$-dimensional metric and an $1$-dimensional metric, the property  
(\ref{4.3.00444}) also follows from 
Corollary 4.1 of \cite{DEF1}. 
As it was stated in Section 3, any semi-Riemannian manifold satisfying (\ref{4.3.00444})
is Riemann compatible. Thus the G\"{o}del metric satisfies also (\ref{A09}). Now (\ref{A09}), by (\ref{aaa01}), turns into 
$\ \omega _{r} g^{rs} ( \omega _{h}R_{sijk} + \omega _{j}R_{sikh} +  \omega _{k}R_{sihj} ) = 0$, 
where $R_{sijk}$ and $g^{rs}$ are the local components of the the Riemann-Christoffel curvature tensor $R$ and the tensor $g^{-1}$ 
of the G\"{o}del metric $g$. We note
that the $1$-form $\omega $, with respect to Definition 3.1 of \cite{2233}, is called $R$-compatible and hence Weyl compatible.
Thus we have

\begin{theorem}\label{thm4.1}
The G\"{o}del spacetime $(M,g)$ is a cyclic Ricci parallel and Riemann compatible 
manifold satisfying:
$\mbox{rank}\, S\, =\, 1$,
$R \cdot R = Q(S,R)$,
$conh(R) \cdot C = conh(R) \cdot conh(R) = 0$,
$C \cdot C = C \cdot conh(R) = \frac{\kappa }{6} \, Q(g, C)$ 
and the $1$-form $\omega $, defined by (\ref{aaa01}), is $R$-compatible as well as Weyl compatible.
\end{theorem}
The above presented results lead to the following generalizations.

\indent Let
$(\overline{M} \times \widetilde{N},g = \overline{g} \times \widetilde{g})$ 
be the product manifold of an $(n-1)$-dimensional semi-Riemannian manifold $(\overline{M},\overline{g})$, $n \geq 4$,
and an $1$-dimensional manifold $(\widetilde{N},\widetilde{g})$.
Moreover, let $(\overline{M},\overline{g})$ be a conformally flat manifold, provided that $n \geq 5$.
The local components $C_{hijk}$, $h,i,j,k \in \{1,2, \ldots , n \}$, 
of the Weyl conformal curvature tensor $C$ of $(\overline{M} \times \widetilde{N},g)$ which may not vanish
identically are the following (cf. \cite{D6}, eqs. (49)-(51)) 
\begin{equation}
\label{4.3.010}
C_{abcd} = \frac{1}{(n-3)(n-2)}\, ( \overline{g}_{ad}A_{bc} - \overline{g}_{ac}A_{bd} 
+ \overline{g}_{bc}A_{ad} -  \overline{g}_{bd}A_{ac}) ,
\end{equation}
\begin{equation}
\label{4.3.011}
C_{nbcn} = - \frac{1}{n-2}\, \widetilde{g}_{nn}A_{bc} ,
\end{equation}
where 
$A_{ab} = \overline{S}_{ab} - \frac{ \overline{\kappa} }{n-1} \, \overline{g}_{ab}$,
and $\overline{g}_{ab}$ and $\overline{S}_{ab}$ denote the local components 
of the metric tensor $\overline{g}$ and the Ricci tensor $\overline{S}$ of $(\overline{M},\overline{g})$, respectively,
$a, b, c, d \in \{ 1, 2, \ldots , n-1 \}$, and  
$\overline{\kappa }$ is the scalar curvature of $(\overline{M},\overline{g})$.
Further, we denote by $U_{C}$ the set of all points 
of $(\overline{M} \times \widetilde{N},g)$ 
at which the Weyl conformal curvature tensor $C$ of $(\overline{M} \times \widetilde{N},g)$ is non-zero.
We note that the tensor $C$ is non-zero at a point of $U_{C}$ 
if and only if $\overline{S} \neq \frac{ \overline{\kappa} }{n-1} \, \overline{g}$ at this point.

As an immediate consequence of Proposition 2 of \cite{D6} we get the following equivalence:
(\ref{4.3.012})
holds on the set $U_{C}$ of the defined above manifold $(\overline{M} \times \widetilde{N},g)$
for some function $L_{C}$ on this set, 
if and only if at every point of $U_{C}$ we have
\begin{equation}
\label{4.3.014}
g^{bc}A_{ab}A_{cd} = (n-3)(n-2) L_{C}\, A_{ad} + \lambda \, \overline{g}_{ad} ,
\end{equation}
\begin{equation}
\label{4.3.015}
B_{ad}B_{bc} - B_{ac}B_{bd} 
= \left( (n-2)^{2}L_{C}^{2} - \frac{\lambda }{n-1}\right)\, (\overline{g}_{ad}\overline{g}_{bc} - \overline{g}_{ac}\overline{g}_{bd}) ,
\end{equation}
where 
$B = A + (n-2)L_{C}\, \overline{g}$,
and  
$\lambda $ is a constant. Furthermore, in view of Lemma 3.1 of \cite{GLOG}, at every point of $U_{C}$
(\ref{4.3.015}) is equivalent to $\mbox{rank}\, B = 1$, i.e. 
\begin{equation}
\label{4.3.016}
\mbox{rank} \left( \overline{S} - \left(\frac{ \overline{\kappa }}{n-1} - (n-2) L_{C}\right)\, \overline{g} \right) = 1 .
\end{equation}

Further, we denote by $S$ and $\kappa$ the Ricci tensor and the scalar curvature of   
$(\overline{M} \times \widetilde{N},g = \overline{g} \times \widetilde{g})$, respectively.
It is obviuos that 
$S = \overline{S}$ and $\kappa = \overline{\kappa }$. Therefore (\ref{4.3.016}) yields
\begin{equation}
\label{4.3.0166}
\mbox{rank} \left(S - \left(\frac{ \kappa }{n-1} - (n-2) L_{C}\right)\,g \right) = 1.
\end{equation}
From the above presented considerations and Proposition 3.1 it follows
\begin{theorem}
Let
$(\overline{M} \times \widetilde{N},g = \overline{g} \times \widetilde{g})$ be
the product manifold of an $(n-1)$-dimensional semi-Riemannian manifold $(\overline{M},\overline{g})$, $n \geq 4$,
and an $1$-dimensional manifold $(\widetilde{N},\widetilde{g})$. Moreover, 
let $(\overline{M},\overline{g})$ be a conformally flat manifold, provided that $n \geq 5$.
If on $\overline{M}$ we have
$\mbox{rank}\, ( \overline{S} - \rho \, \overline{g} ) = 1$, for some function $\rho$, then 
$\mbox{rank}\, ( S - \rho \, g ) = 1$ and (\ref{4.3.012}), i.e.
$C \cdot C = L_{C}\, Q(g,C)$,
with $L_{C} = \frac{1}{n-2} ( \frac{ \kappa }{n-1} - \rho )$,
hold on $\overline{M} \times \widetilde{N}$.
In particular, if the rank of the Ricci tensor of $(\overline{M},\overline{g})$ is one,
then the rank of the Ricci tensor of
$\overline{M} \times \widetilde{N}$
is also one
and (\ref{4.3.012}), with $L_{C} = \frac{ \kappa }{(n-2)(n-1)}$,
or equivalently,
$conh(R) \cdot conh(R) = 0$
holds on this manifold. 
\end{theorem}

We present now an application of the last theorem.

From Theorem 4.1 of \cite{32} it follows that a hypersurface $M$ 
immersed isometrically in a semi-Riemannian space of constant curvature $N$, $\dim N \geq 5$,  
is a quasi-umbilical hypersurface if and only if it is a conformally flat manifold.
Furthermore, using the Gauss equation of $M$ in $N$, we can easily prove that if $M$ is quasi-umbilical hypersurface
then it is also a quasi-Einstein manifold. These facts, together with Theorem 4.2, leads to the following

\begin{theorem}
Let $(\overline{M},\overline{g})$ be a manifold
which is isometric with a quasi-umbilical hypersurface immersed isometrically 
in a semi-Riemannian space of constant curvature $N$, $\dim N \geq 5$. 
Let $(\widetilde{N},\widetilde{g})$ be an $1$-dimensional manifold. 
Then the manifold
$(\overline{M} \times \widetilde{N},g = \overline{g} \times \widetilde{g})$ 
is a quasi-Einstein manifold with pseudosymmetric Weyl conformal curvature tensor.
\end{theorem}

In this way we obtain a family of quasi-Einstein manifolds 
with pseudosymmetric Weyl conformal curvature tensor.
We mention that quasi-Einstein hypersurfaces 
with pseudosymmetric Weyl conformal curvature tensor
immersed isometrically 
in semi-Riemannian spaces of constant curvature
were investigated in \cite{GLOG2}.

\indent Theorem 4.3 together with Theorem 4.2 and Example 4.1 of \cite{49}, yields
\begin{theorem}\label{thm4.4}
Let $(\overline{M},\overline{g})$ be a manifold
which is isometric with a quasi-umbilical hypersurface immersed isometrically 
in a semi-Euclidean space $N$, $\dim N \geq 5$. 
Let $(\widetilde{N},\widetilde{g})$ be an $1$-dimensional manifold. 
Then the manifold
$(\overline{M} \times \widetilde{N},g = \overline{g} \times \widetilde{g})$ 
is a quasi-Einstein manifold with pseudosymmetric Weyl conformal curvature tensor satisfying $R \cdot R = Q(S,R)$.
\end{theorem}

\indent Finally, we consider some extension of the G\"{o}del metric.
\newline
{\bf{Example 4.1.}} 
(i) We define the metric $g$ on $M \, =\, \{ (t,r,\phi ,z): t > 0, r > 0 \} \subset {\mathbb R}^{4}$ by 
(cf. \cite{ReboTio}, Section 1)
\begin{eqnarray}
ds^{2} 
&=& 
(dt + H(r)\, d\phi )^{2}
- D^{2}(r)\, d\phi ^{2}
- dr^{2}
- dz^{2} ,
\label{godel77a}
\end{eqnarray}
where $H$ and $D$ are certain functions on $M$. 
In the special case, if 
$H(r) \, =\, \frac{2 \sqrt{2}}{m} \, sinh^{2} (\frac{m r}{2})$ and 
$D(r)\, =  \, \frac{2}{m} \, sinh (\frac{m r}{2})\, cosh (\frac{m r}{2})$ 
then $g$ is the G\"{o}del metric (e.g. see \cite{ReboTio}, eq. (1.6)).
\newline
(ii) Since the metric
$g$ 
defined by
(\ref{godel77a})
is the product metric of a $3$-dimensional metric 
and a $1$-dimensional metric
(\ref{4.3.00444})
holds on $M$.
We can check that the Riemann-Christoffel curvature tensor $R$ of $(M,g)$ is expressed by 
a linear combination of 
the Kulkarni-Nomizu products formed by $S$ and $S^{2}$, i.e. by the tensors 
$S \wedge S$, $S \wedge S^{2}$ and $S^{2} \wedge S^{2}$,
\begin{eqnarray*}
R &=& \phi _{1} \, S \wedge S + \phi _{2} \, S \wedge S^{2} + \phi _{3} \, S^{2} \wedge S^{2} ,
\end{eqnarray*}
\begin{eqnarray*}
\phi _{1} 
&=& 
\frac{D^{2}}{\tau } ( 2 D^{2} H''^{2} - 4 D D' H' H'' - 3 H'^{4} 
+ 8 D D'' H'^{2} + 2 D'^{2} H'^{2} - 8 D^{2} D''^{2} ), \\
\phi _{2} 
&=& 
\frac{2 D^{4}}{\tau} ( H'^{2} - 4 D D'' ),\\
\phi _{3} &=&
 - \frac{4 D^{6}}{\tau }, \ \ H' \ =\ \frac{d H}{d r}, \ \ H'' \ =\ \frac{d H' }{d r} ,
\\
\tau  
&=& 
( H'^{2} - 2 D D'' ) 
( D^{2} H''^{2} - 2 D D' H' H'' - H'^{4} + 2 D D'' H'^{2} + D'^{2} H'^{2} ) ,
\end{eqnarray*}
provided that the function $\tau $  
is non-zero at every point of $M$. 
\newline
(iii) If $H(r)\, =\, a r^{2}$, $a\, =\, const. \neq 0$ and $D(r)\, =\, r$  
then (\ref{godel77a}) turns into (\cite{ReboTio}, eq. (3.20))
\begin{eqnarray}
ds^{2} & = & 
(dt +  a r^{2}\, d\phi )^{2}
- r^{2}\, d\phi ^{2}
- dr^{2}
- dz^{2} .
\label{newgodel77a}
\end{eqnarray}
The spacetime $(M,g)$ with the metric $g$ defined by (\ref{newgodel77a})
is called the Som-Raychaudhuri solution of the Einstein field
equations (\cite{SomRay}). For the metric (\ref{newgodel77a}) the function $\tau$ is non-zero
at every point of $M$. 
\newline
(iv) 
We refer to \cite{DGHSaw} and \cite{Malgorz1} for surveys on 
semi-Riemannian manifolds
$(M,g)$, $n \geq 4$, having Riemann-Christoffel curvature tensor $R$ expressed by a linear combination of 
the Kulkarni-Nomizu products formed by $g$ and $S$, i.e. by the tensors 
$g \wedge g$, $g \wedge S$ and $S \wedge S$. In particular, we mention that 
in the class of the Reissner-Nordstr\"{o}m-de Sitter spacetimes there are spacetimes having that property
(\cite{Kow2}, Example 5.3).
\newline

\indent
It may be mentioned that we have calculated the local components of various tensors using Wolfram Mathematica,
as well as 
SymPy and Maxima packages for symbolic calculation.
\section*{Conclusion:}
By considering the dust particles as galaxies, the G\"{o}del spacetime can be taken 
as a cosmological model of rotating universe. Although G\"{o}del spacetime is not a realistic model of the universe 
in which we live but it realized many peculiar properties. For example, the existence of closed timelike curves implies 
a form of time travel in an alternative universe described by the G\"{o}del spacetime. Also G\"{o}del spacetime 
is quasi-Einstein, Ricci tensor is cyclic parallel but not Codazzi type, which may be physically interpreted as the content 
of the spacetime is of rotating matter without singularity. It is neither pseudosymmetric 
nor Ricci pseudosymmetric but a special type of Ricci generalized pseudosymmetric, 
and it is not conformally pseudosymmetric but its Weyl conformal curvature tensor is pseudosymmetric 
(i.e., $C \cdot C = \frac{\kappa}{6}Q(g,C)$) and also the spacetime is Riemann compatible as well as Weyl 
compatible (Theorem \ref{thm4.1}). Hence G\"{o}del spacetime forced us 
to obtain a new class of semi-Riemannian manifolds which is quasi-Einstein with pseudosymmetric Weyl conformal curvature tensor 
and is a special type of Ricci generalized pseudosymmetric manifolds (Theorem \ref{thm4.4}).

\section*{Acknowledgments}
The first and third named authors are supported by a grant of the Wroc\l aw University of Environmental 
and Life Sciences, Poland [WIKSiG/441/212/S]. The fourth and fifth named authors gratefully acknowledge the financial support of CSIR, New Delhi, India [File no: 09/025(0194)/2010-EMR-I, Project F. No. 25(0171)/09/EMR-II].

\end{document}